\documentclass[12pt]{article}
\usepackage{color}
\usepackage{graphicx}
\usepackage{colortbl}
\usepackage{array}
\usepackage{amssymb}
\usepackage{amsmath}
\usepackage{amsthm}
\usepackage{mathrsfs}
\usepackage{dcolumn}
\usepackage{longtable}
\usepackage{hhline}
\usepackage{cancel}
\usepackage{fancyhdr}
\usepackage{mathtools}
\numberwithin{equation}{section}
\allowdisplaybreaks

\newtheorem{thm}{Theorem}[section]
\newtheorem{defn}[thm]{Definition}

\title{\bf The Multivariable Generalized Hermite-Type-Genocchi Polynomials of Order $a$}
\author{\textbf{Roberto B. Corcino$^{1,2}$}\\{\large\bf Cristina B. Corcino$^{1,2}$}\\{$^{1}$Research Institute for Computational}\\ {Mathematics and Physics}
	\\$^{2}$Mathematics Department\\Cebu Normal University\\Cebu City, Philippines \vspace{9pt}}

\date{}

\begin{document}
	
	\maketitle
	
		\begin{abstract}	
		This study presents a new class of poly-Genocchi polynomials constructed through the integration of some interesting polynomials. The resulting family, referred to as the multivariable generalized Hermite-type-Genocchi polynomials of order $a$, is investigated in detail. Several fundamental properties are derived, including explicit representations, addition formulas, and polynomial expansions. In addition, relationships between this new family of polynomials and a certain generalized Stirling numbers of the first kind are established.
		
		
				
	\end{abstract} 
	
	\section{Introduction}
	The multivariable Hermite polynomials, denoted by $\mathcal{D}^{[m]}_n(h_1,h_2,\ldots,h_m)$, is defined \cite{Dattoli-0} as follows:
	\begin{equation}\label{multi_hermite}
		\exp(h_1t+h_2t^2+\ldots+h_mt^m)=\sum_{l=0}^{\infty}\mathcal{D}^{[m]}_l(h_1,h_2,\ldots,h_m)\frac{t^l}{l!}.
	\end{equation}
	When $m=2$, equation \eqref{multi_hermite} gives the following extended Hermite polynomials \cite{Dattoli-1,Khan}:
	\begin{equation*}
		e^{pt^2+qt}=\sum_{n=0}^{\infty}H_n(p,q)\frac{t^n}{n!}.
	\end{equation*}
	With $q=-1$ and $p=2x$, we obtain the well-known classical Hermite polynomials:
	\begin{equation*}
		\sum_{l=0}^{\infty}H_l(x)\frac{t^l}{l!}=e^{2xt^2-t}.
	\end{equation*}
	When $m=3$, equation \eqref{multi_hermite} gives the three-variable Hermite polynomials \cite{ref2-1,ref2-2}. With this latter generalization, Araci et al. \cite{ref2-1} introduced the following variation of Genocchi polynomials as follows:
	\begin{equation*}
		 \sum_{l=0}^{\infty}{}_H\mathcal{G}^{(a)}_l(q,p,r;\lambda,u)\frac{t^l}{l!}=e^{qt+pt^2+rt^3}\left(\frac{(1-u)t}{\lambda e^{t}-u}\right)^{a}.
	\end{equation*}
	On the other hand, the degenerate Hermite polynomials were introduced by T. Kim \cite{ref45-2} as follows:
	$$\sum_{l=0}^{\infty}H_{l,\lambda}(p)\frac{t^l}{l!}=e^{-1}_{\lambda}(t)e^{p}_{\lambda}(2t)$$
	where $e^{p}_{\lambda}(t)$ equals
	$$e^{p}_{\lambda}(x)=(\lambda x+1)^{\frac{p}{\lambda}}= \sum_{l=0}^{\infty}(p)_{l,\lambda}\frac{x^l}{l!},$$
	which satisfy the following properties
	\begin{align*}
		e^{p+q}_{\lambda}(x)&=e^p_{\lambda}(x)e^q_{\lambda}(x),\\
		\frac{d}{dz}e^z_{\lambda}(x)&=\log (\lambda x+1)^{1/\lambda}e^z_{\lambda}(x),
	\end{align*}
	where $(p)_{n,\lambda}$ denotes the degenerate falling factorial. 
	
	\smallskip
	The degenerate version of Genocchi, Euler and Bernoulli polynomials were introduced by Carlitz and Lim \cite{Carlitz-2, Lim}. Kim-Kim \cite{Kim-4, ref45-1-3} introduced the following:
	\begin{align*}
		\frac{1-u}{e_{\lambda}(x)-u}e^z_{\lambda}(x)&=\sum_{l=0}^{\infty}h_{l,\lambda}(z|u)\frac{x^l}{l!}\\
		\frac{2x^r}{e_{\lambda}(x)+1}e^z_{\lambda}(x)&=\sum_{l=0}^{\infty}A^{(r)}_{l,\lambda}(z)\frac{x^l}{l!},
	\end{align*}
	by integrating the Genocchi, degenerate exponential and Euler polynomials. More properties are obtained related to  $h_{l,\lambda}(x|u)$ \cite{ref45-1-4}.
		
	\smallskip
	The degenerate version of poly-Euler polynomials were introduced \cite{ref51-1} as follows:
	\begin{equation*}
		2\frac{{\rm Ei}_{k,\nu}(\log_{\nu}(1+x))}{xe_{\nu}(x)+1}e^z_{\nu}(x)=\sum_{l=0}^{\infty}\mathcal{E}^{(k)}_{l,\nu}(z)\frac{x^l}{l!},
	\end{equation*}
	where $k \in \mathbb{Z}$ and
	\begin{align}
		&{\rm Ei}_{k,\nu}(z)=\sum_{l=1}^{\infty}\frac{(1)_{l,\nu}z^l}{l^k(l-1)!}, \;\;\; \nu\in \mathbb{R},\label{poly_exp-1}\\
		&\frac{d}{dz}{\rm Ei}_{k,\nu}(\log_{\nu}(1+z))=\frac{(1+z)^{\nu-1}}{\log_{\nu}(1+z)}{\rm Ei}_{k-1,\nu}(\log_{\nu}(1+z)),\nonumber\\
		&{\rm Ei}_{k,\nu}(\log_{\nu}(1+z))\nonumber\\
		&\qquad =\int_0^z\frac{(1+z)^{\nu-1}}{\log_{\nu}(1+z)}\int_0^z\ldots \frac{(1+z)^{\nu-1}}{\log_{\nu}(1+z)}\int_0^z\frac{(1+z)^{\nu-1}}{\log_{\nu}(1+z)}zdz\ldots dz\nonumber\\
		&\qquad=\sum_{m=0}^{\infty}\sum_{s_1+s_2+\ldots +s_{k-1}=m}\binom{m}{s_1,\ldots,s_{k-1}}\nonumber\\
		&\qquad\;\;\;\;\;\times\frac{b_{s_1,\nu}(\nu-1)}{s_1+1}\frac{b_{s_2,\nu}(\nu-1)}{s_1+s_2+1}\ldots \frac{b_{s_{k-1},\nu}(\nu-1)}{s_1+\ldots+s_{k-1}+1}\frac{x^{m+1}}{m!}\label{poly_exp-1-1}
	\end{align}
	the modified degenerate polyexponential function in (\cite{ref45-1, Kim-1, ref45-1-1, ref45-1-2}).
	
	\bigskip
	Using a more general form of the degenerate Hermite polynomials, Corcino-Corcino \cite{Corcino-22} defined the following polynomials:
	\begin{equation}\label{eq711-A1}
		\sum_{n=0}^{\infty}\widehat{\mathcal{G}}_{n}^{(k,a)}(w,z; \nu, \lambda, u, c, d)\frac{x^n}{n!}=\left(\frac{{\rm Ei}_{k,\nu}(\log_{\nu}(1+(1-u)x\ln cd))}{\lambda d^x-uc^{-x}}\right)^{a}e_{\nu}^{w}(x)e_{\nu}^{z}(x^2),
	\end{equation}
	where $|t| < \frac{\sqrt{(\ln \left(\frac{\lambda}{u}\right))^2+4\pi^2}}{|\ln c + \ln c|}.$ 
	When $a=1$, equation \eqref{eq711-A1} yields
	\begin{equation*}
		\sum_{n=0}^{\infty}\widehat{\mathcal{G}}_{n}^{(k)}(w,z; \nu, \lambda, u, c, d)\frac{x^n}{n!}=\frac{{\rm Ei}_{k,\nu}(\log_{\nu}(1+(1-u)x\ln cd))}{\lambda d^t-uc^{-t}}e_{\nu}^{w}(x)e_{\nu}^{z}(x^2).
	\end{equation*}
	where $\widehat{\mathcal{G}}_{n}^{(k)}(w,z; \nu, \lambda, u, c, d)=\widehat{\mathcal{G}}_{n}^{(k,1)}(w,z; \nu, \lambda, u, c, d)$.
	
	\smallskip
	The degenerate three-variable version of these polynomials are defined by Corcino-Corcino \cite{Corcino-23} as follows:
	\begin{align}
		&\sum_{n=0}^{\infty}\widehat{\mathcal{G}}_{n}^{(k,a)}(v,w,z; \nu, \lambda, u, c, d)\frac{x^n}{n!}\nonumber\\
		&\qquad=\left(\frac{{\rm Ei}_{k,\nu}(\log_{\nu}(1+(1-u)x\ln cd))}{\lambda d^t-uc^{-x}}\right)^{a}e_{\nu}^{v}(x)e_{\nu}^{w}(x^2)e_{\nu}^{z}(x^3),\label{eq710-A}
	\end{align}
	where $|t| < \frac{\sqrt{(\ln \left(\frac{\lambda}{u}\right))^2+4\pi^2}}{|\ln c + \ln d|}.$ 
	When $a=1$, equation \eqref{eq710-A} yields
	\begin{align}
		&\sum_{n=0}^{\infty}\widehat{\mathcal{G}}_{n}^{(k)}(v,w,z; \nu, \lambda, u, c, d)\frac{x^n}{n!}\nonumber\\
		&\qquad=\frac{{\rm Ei}_{k,\nu}(\log_{\nu}(1+(1-u)x\ln cd))}{\lambda d^x-uc^{-x}}e_{\nu}^{v}(x)e_{\nu}^{w}(x^2)e_{\nu}^{z}(x^3).\label{eq710-1}
	\end{align}
	where $\widehat{\mathcal{G}}_{n}^{(k)}(w,z; \nu, \lambda, u, c, d)=\widehat{\mathcal{G}}_{n}^{(k,1)}(w,z; \nu, \lambda, u, c, d)$.
	
	\smallskip
	To further generalize these degenerate polynomials, this paper introduces a multivariable extension of the aforementioned degenerate polynomials. This extension is constructed by synthesizing ideas from generalized Genocchi polynomials and multivariable Hermite polynomials.
	
	\smallskip
	In addition, this paper identifies notable identities containing the proposed polynomials and several families of degenerate special numbers \cite{ref-kim-1, ref-kim-2, ref-kim-4, ref-kim-7}. These findings further reflect the expanding scope and active research developments in this domain, as highlighted \cite{ref-kim-5, ref-kim-6, ref-kim-8}.

	\section{Preliminary Concept and Results}

	\begin{defn}\rm The multivariable generalized Hermite-type-Genocchi polynomials of order $a$, denoted by $\widehat{\mathcal{G}}_{n,p}^{(k,a)}(x_1, x_2,\ldots, x_p; \nu, \lambda, u, c, d)$, are defined by
		\begin{align}
			&\sum_{n=0}^{\infty}\widehat{\mathcal{G}}_{n,p}^{(k,a)}(x_1, x_2,\ldots, x_p; \nu, \lambda, u, c, d)\frac{t^n}{n!}\nonumber\\
			&\qquad=\left(\frac{{\rm Ei}_{k,\nu}(\log_{\nu}(1+(1-u)t\ln cd))}{\lambda d^t-uc^{-t}}\right)^{a}\prod_{j=1}^{p}e_{\nu}^{x_j}(t^j),\label{eq711}
		\end{align}
		where $|t| < \frac{\sqrt{(\ln \left(\frac{\lambda}{u}\right))^2+4\pi^2}}{|\ln c + \ln d|}.$ 
		When $a=1$, equation \eqref{eq711} yields
		\begin{align*}
			&\sum_{n=0}^{\infty}\widehat{\mathcal{G}}_{n,p}^{(k)}(x_1, x_2,\ldots, x_p; \nu, \lambda, u, c, d)\frac{t^n}{n!}\\
			&\qquad=\frac{{\rm Ei}_{k,\nu}(\log_{\nu}(1+(1-u)t\ln cd))}{\lambda d^t-uc^{-t}}\prod_{j=1}^{p}e_{\nu}^{x_j}(t^j).
		\end{align*}
		where $\widehat{\mathcal{G}}_{n,p}^{(k)}(x_1, x_2,\ldots, x_p; \nu, \lambda, u, c, d)=\widehat{\mathcal{G}}_{n}^{(k,1)}(x_1, x_2,\ldots, x_p; \nu, \lambda, u, c, d)$ denotes the Hermite-type multivariable generalized poly-Genocchi polynomials.
	\end{defn}
	
	\bigskip
	\begin{thm}\label{thm-2-2} The multivariable generalized Hermite-type-Genocchi polynomials are equal to
		\begin{align*}
			&\frac{1}{l+1}\widehat{\mathcal{G}}_{l+1,p}^{(k)}(x_1, x_2,\ldots, x_p; \nu, \lambda, u, c, d)\\
			&\;\;\;\;\;=\sum_{q=0}^{l}\sum_{i=0}^{l-q}\sum_{j=0}^{l-q-i}\binom{l}{q}\binom{l-q}{i}\binom{l-q-i}{j}
			\left(\sum_{n_1+n_2+\ldots+n_p= i}\,i!\prod_{j=1}^{p}\nu^{\frac{n_j}{j}} \binom{\frac{x_j}{\nu}}{\frac{n_j}{j}}\right)\\
			&\qquad\times\left(\frac{-1}{u}\right) \,\,\frac{(\log c)^{j}\,(\log cd)^{l-q-i-j}}{\left(1-\frac{\lambda}{u}\right)^{l-q-i-j+1}}\mathbb{B}_{s_1,s_2,\ldots,s_{k-1}}(q,\nu-1)A_{l-q-i-j}\left(\frac{u}{\lambda}\right).
		\end{align*}
		where $\widehat{\mathcal{G}}_{0}^{(k,a)}(x_1, x_2,\ldots, x_p; \nu, \lambda, u, c, d)=0$, $A_m(s)$ denotes the Eulerian polynomial and $\mathbb{B}_{s_1,s_2,\ldots,s_{k-1}}(q,\nu-1)$ is defined in equation \eqref{deg_bern-multi-1}.
		\begin{proof}
			Using equation \eqref{poly_exp-1-1} with $x=(1-u)t\ln cd$ and the following identities and notation,
			\begin{align}
				&e^{x_p}_\nu (t^p)=\sum_{m=0\atop p|m}^{\infty}\nu^{\frac{m}{p}} m!\binom{\frac{x_p}{\nu}}{\frac{m}{p}}\, \frac{t^m}{m!}\,\,\,\,\,\,,\binom{\frac{x_p}{\nu}}{i}=0\,\,,i \notin \mathbb{N }\cup\{0\},\nonumber\\
				&\frac{\prod_{j=1}^{p}e_{\nu}^{x_j}(t^j)}{\lambda d^t -uc^{-t}}=\sum_{l=0}^{\infty}\left\{\sum_{i=0}^{l}\binom{l}{i}\left(\sum_{n_1+n_2+\ldots+n_p= i}\,i!\prod_{j=1}^{p}\nu^{\frac{n_j}{j}} \binom{\frac{x_j}{\nu}}{\frac{n_j}{j}}\right)\right.\nonumber\\
				&\;\;\;\;\;\left.\times\left(-\frac{1}{u}\sum_{j=0}^{n-i}\,\,\sum_{s=0}^{\infty}\binom{l-i}{j}\left(\frac{\lambda}{u}\right)^s \, \,(\log c)^{j}\,(s\log cd)^{n-i-j}\right)\right\}\frac{t^l}{l!}\nonumber\\
				&\mathbb{B}_{s_1,s_2,\ldots,s_{k-1}}(m,\nu-1)\nonumber\\
				&\;\;\;\;\;=((1-u)\ln cd)^{m+1}\sum_{s_1+s_2+\ldots +s_{k-1}=m}\binom{m}{s_1,\ldots,s_{k-1}}\nonumber\\
				&\;\;\;\;\;\;\;\;\;\;\times\frac{b_{s_1,\nu}(\nu-1)}{s_1+1}\frac{b_{s_2,\nu}(\nu-1)}{s_1+s_2+1}\ldots \frac{b_{s_{k-1},\nu}(\nu-1)}{s_1+\ldots+s_{k-1}+1},\label{deg_bern-multi-1}
			\end{align}
		we have
		\begin{align*}
			&\sum_{l=1}^{\infty}\frac{1}{l}\widehat{\mathcal{G}}_{l,p}^{(k)}(x_1, x_2,\ldots, x_p; \nu, \lambda, u, c, d)\frac{t^{l-1}}{(l-1)!}\\
			&\qquad=\sum_{l=0}^{\infty}\frac{1}{l+1}\widehat{\mathcal{G}}_{l+1,p}^{(k)}(x_1, x_2,\ldots, x_p; \nu, \lambda, u, c, d) \frac{t^{l}}{l!}\\
			&\qquad=\sum_{l=0}^{\infty}\left\{\sum_{q=0}^{l}\sum_{i=0}^{l-q}\sum_{j=0}^{l-q-i}\,\,\sum_{s=0}^{\infty}\binom{l}{q}\binom{l-q}{i}\binom{l-q-i}{j}\right.\\
			&\left.\qquad\qquad\times\left(\sum_{n_1+n_2+\ldots+n_p= i}\,i!\prod_{j=1}^{p}\nu^{\frac{n_j}{j}} \binom{\frac{x_j}{\nu}}{\frac{n_j}{j}}\right)\left(\frac{-1}{u}\right)\left(\frac{\lambda}{u}\right)^s \,\,(\log c)^{j}\right.\\
			&\left.\qquad\qquad\times(s\log cd)^{l-q-i-j}\mathbb{B}_{s_1,s_2,\ldots,s_{k-1}}(q,\nu-1)\right\}\frac{t^l}{l!}.
		\end{align*}
		Thus,   
		\begin{align*}
			&\frac{1}{l+1}\widehat{\mathcal{G}}_{l+1,p}^{(k)}(x_1, x_2,\ldots, x_p; \nu, \lambda, u, c, d)\\
			&\;\;\;\;\;=\sum_{q=0}^{n}\sum_{i=0}^{l-q}\sum_{j=0}^{l-q-i}\binom{l}{q}\binom{l-q}{i}\binom{l-q-i}{j}
			\left(\sum_{n_1+n_2+\ldots+n_p= i}\,i!\prod_{j=1}^{p}\nu^{\frac{n_j}{j}} \binom{\frac{x_j}{\nu}}{\frac{n_j}{j}}\right)\\
			&\qquad\times\left(\frac{-1}{u}\right) \,\,(\log c)^{j}\,(\log cd)^{l-q-i-j}\mathbb{B}_{s_1,s_2,\ldots,s_{k-1}}(q,\nu-1)\frac{A_{l-q-i-j}\left(\frac{u}{\lambda}\right)}{\left(1-\frac{\lambda}{u}\right)^{l-q-i-j+1}},
		\end{align*}
		where $$\sum_{s=0}^{\infty}\left(\frac{\lambda}{u}\right)^ss^{l-q-i-j}=\frac{A_{l-q-i-j}\left(\frac{u}{\lambda}\right)}{\left(1-\frac{\lambda}{u}\right)^{l-q-i-j+1}},$$
		(see \cite[p.245]{ref22}), with $A_l(u)$ denotes the Eulerian polynomials.
		\end{proof}
	\end{thm}
	
	\bigskip
	\begin{thm}\label{thm-2-3} The multivariable generalized Hermite-type-Genocchi polynomials of order $a$ are equal to
		\begin{align*}
			&\frac{1}{(n+a)_{a}}\widehat{\mathcal{G}}_{n+a,p}^{(k,a)}(x_1, x_2,\ldots, x_p; \nu, \lambda, u, c, d)\\
			&\;\;\;\;\;\;\;\;\;\;=\sum_{q=0}^{n}\binom{n}{q}
			\left(\sum_{n_1+n_2+\ldots+n_p= n-q}\,(n-q)!\prod_{j=1}^{p}\nu^{\frac{n_j}{j}} \binom{\frac{x_j}{\nu}}{\frac{n_j}{j}}\right)\\
			&\;\;\;\;\;\;\;\;\;\;\;\;\;\;\;\times\sum_{k_1+k_2+\ldots+k_{a}=n-q}\prod_{i=1}^{a}\left\{\sum_{j=0}^{k_i}\sum_{p=0}^{\infty}\frac{1}{k_i!}\binom{k_i}{j}\right.\\
			&\left.\;\;\;\;\;\;\;\;\;\;\;\;\;\;\;\times\,\mathbb{B}_{s_1,s_2,\ldots,s_{k-1}}(j,\nu-1)\left(\frac{u}{\lambda}\right)^p(-p\log cd)^{k_i-j}\right\}.
		\end{align*}
		where $\widehat{\mathcal{G}}_{n}^{(k,a)}(x_1, x_2,\ldots, x_p; \nu, \lambda, u, c, d)=0$ for $n=0, 1, \ldots, a-1$, $A_l(u)$ denotes the Eulerian polynomial and $\mathbb{B}_{s_1,s_2,\ldots,s_{k-1}}(i,\nu-1)$ is defined in equation \eqref{deg_bern-multi-1}.
		\begin{proof}
			Raising to power $a$, the following quotient,
			\begin{align*}
				&\frac{{\rm Ei}_{k,\nu}(\log_{\nu}(1+(1-u)t\ln cd))}{\lambda d^t-uc^{-t}}\\
				&\;\;\;\;\;=t\sum_{m=0}^{\infty}\left\{\sum_{j=0}^{m}\sum_{n=0}^{\infty}\frac{1}{m!}\binom{m}{j}\mathbb{B}_{s_1,s_2,\ldots,s_{k-1}}(j,\nu-1)\left(\frac{u}{\lambda}\right)^n(-n\log cd)^{m-j}\right\}t^m,
			\end{align*}
			gives
			\begin{align*}
				&\left(\frac{{\rm Ei}_{k,\nu}(\log_{\nu}(1+(1-u)t\ln cd))}{\lambda d^t-uc^{-t}}\right)^{a}\\
				&\;\;\;\;\;=t^{a}\sum_{m=0}^{\infty}\sum_{k_1+k_2+\ldots+k_{a}=m}\prod_{i=1}^{a}\left\{\sum_{j=0}^{k_i}\sum_{n=0}^{\infty}\frac{1}{k_i!}\binom{k_i}{j}\mathbb{B}_{s_1,s_2,\ldots,s_{k-1}}(j,\nu-1)\right.\\
				&\;\;\;\;\;\;\;\;\;\;\left.\left(\frac{u}{\lambda}\right)^n(-n\log cd)^{k_i-j}\right\}t^m.
			\end{align*}
		Thus,
			\begin{align*}
				&\sum_{n=-a}^{\infty}\frac{1}{(n+a)_{a}}\widehat{\mathcal{G}}_{n+a,p}^{(k,a)}(x_1, x_2,\ldots, x_p; \nu, \lambda, u, c, d)\frac{t^{n}}{n!}\\
				&\;\;\;\;\;\;\;\;\;\;=\sum_{n=0}^{\infty}\sum_{q=0}^{n}\binom{n}{q}
				\left(\sum_{n_1+n_2+\ldots+n_p= n-q}\,(n-q)!\prod_{j=1}^{p}\nu^{\frac{n_j}{j}} \binom{\frac{x_j}{\nu}}{\frac{n_j}{j}}\right)\\
				&\;\;\;\;\;\;\;\;\;\;\;\;\;\;\;\times\sum_{k_1+k_2+\ldots+k_{a}=n-q}\prod_{i=1}^{a}\left\{\sum_{j=0}^{k_i}\sum_{p=0}^{\infty}\frac{1}{k_i!}\binom{k_i}{j}\right.\\
				&\left.\;\;\;\;\;\;\;\;\;\;\;\;\;\;\;\times\,\mathbb{B}_{s_1,s_2,\ldots,s_{k-1}}(j,\nu-1)\left(\frac{u}{\lambda}\right)^p(-p\log cd)^{k_i-j}\right\}\frac{t^n}{n!},
			\end{align*}
		which consequently proves the theorem.
		\end{proof}
	\end{thm}
		
	\section{Further Results}
	
	Assigning specific values to some parameters yields the following interesting polynomials: 
	\begin{enumerate}
		\item Using equation \eqref{poly_exp-1}, when $k=1$, equation \eqref{eq711} yields
		\begin{equation*}
			\sum_{n=0}^{\infty}\widehat{\mathcal{G}}_{n,p}^{(a)}(x_1, x_2,\ldots, x_p; \nu, \lambda, u, c, d) \frac{t^n}{n!}=\left(\frac{(1-u)t\ln cd}{\lambda d^t-uc^{-t}}\right)^{a}\prod_{j=1}^{p}e_{\nu}^{x_j}(t^j),
		\end{equation*}
		Moreover,
		\begin{equation*}
			\sum_{n=0}^{\infty}\widehat{\mathcal{G}}_{n,p}^{(1)}(x_1, x_2,\ldots, x_p; \nu, \lambda, u, c, d) \frac{t^n}{n!}=\frac{(1-u)t\ln cd}{\lambda d^t-uc^{-t}}\prod_{j=1}^{p}e_{\nu}^{x_j}(t^j).
		\end{equation*} 
		\item For $x_i=0,\,\forall i=1,\ldots, p$, we have
		\begin{equation}\label{eq13-1}
			\sum_{n=0}^{\infty}\widehat{\mathcal{G}}_{n}^{(k,a)}(\nu, \lambda, u, c, d)\frac{t^n}{n!}=\left(\frac{{\rm Ei}_{k,\nu}(\log_{\nu}(1+(1-u)t\ln cd))}{\lambda d^t-uc^{-t}}\right)^{a},
		\end{equation}
		\item For $c=1,d=e$, equation \eqref{eq711} will reduce to
		\begin{equation}
			\sum_{n=0}^{\infty}\widehat{\mathcal{G}}_{n,p}^{(k,a)}(x_1, x_2,\ldots, x_p; \nu, \lambda, u) \frac{t^n}{n!}=\left(\frac{{\rm Ei}_{k,\nu}(\log_{\nu}(1+(1-u)t))}{\lambda e^{t}-u}\right)^{a} \prod_{j=1}^{p}e_{\nu}^{x_j}(t^j).\label{ApostolGenoPoly-1}
		\end{equation}
		When $x_i=0,\,\forall i=1,\ldots, p$, we get
		\begin{equation*}
			\sum_{n=0}^{\infty}\widehat{\mathcal{G}}_{n}^{(k,a)}(\lambda, \nu, u) \frac{t^n}{n!}=\left(\frac{{\rm Ei}_{k,\nu}(\log_{\nu}(1+(1-u)t)}{\lambda e^{t}-u}\right)^{a}.
		\end{equation*}
		\item As $\nu\to 0$, equation \eqref{eq711} approaches to
		\begin{align}
			&\sum_{n=0}^{\infty}\widehat{\mathcal{G}}_{n,p}^{(k,a)}(x_1, x_2,\ldots, x_p;  0, \lambda, u, c, d) \frac{t^n}{n!}\nonumber\\
			&\qquad=\left(\frac{{\rm Ei}_{k,0}(\log_{0}(1+(1-u)t\ln cd))}{\lambda b^t-ua^{-t}}\right)^{a}\prod_{j=1}^{p}e_{0}^{x_j}(t^j)\nonumber\\
			&\sum_{n=0}^{\infty}\widehat{\mathcal{G}}_{n,p}^{(k,a)}(x_1, x_2,\ldots, x_p; \nu, \lambda, u, c, d) \frac{t^n}{n!}\nonumber\\
			&\qquad=\left(\frac{{\rm Ei}_{k}(\log(1+(1-u)t\ln cd))}{\lambda d^t-uc^{-t}}\right)^{a}e^{\sum_{j=1}^{p}x_jt^j},\label{eq13}
		\end{align}
		\item For $\lambda=1$, equation \eqref{ApostolGenoPoly-1} yields 
		\begin{align}
			&\sum_{n=0}^{\infty}\widehat{\mathcal{G}}_{n,p}^{(k,a)}(x_1, x_2,\ldots, x_p; u, 1, e) \frac{t^n}{n!}\nonumber\\
			&\qquad=\left(\frac{{\rm Ei}_{k,\nu}(\log_{\nu}(1+(1-u)t))}{e^{t}-u}\right)^{a} \prod_{j=1}^{p}e_{\nu}^{x_j}(t^j).\label{ApostolGenPoly}
		\end{align} 
		\item For $k=1,\, p=3$, $c=1$ and $d=e$, equation \eqref{eq13} gives 
		\begin{equation}
			\sum_{n=0}^{\infty}\widehat{\mathcal{G}}_{n}^{(1,a)}(v,w,z; \lambda, u) \frac{t^n}{n!}=\left(\frac{(1-u)t}{\lambda e^{t}-u}\right)^{a} e^{vt+wt^2+zt^3},\label{ApostolGenPoly1}
		\end{equation}
		and when $\lambda=1$, equation \eqref{ApostolGenPoly1} gives 
		\begin{equation*}
			\sum_{n=0}^{\infty}\widehat{\mathcal{G}}_{n}^{(1,a)}(v,w,z; 1, u) \frac{t^n}{n!}=\left(\frac{(1-u)t}{e^{t}-u}\right)^{a} e^{vt+wt^2+zt^3},
		\end{equation*}
		Furthermore, when $a=1$, we have
		\begin{equation*}
			\sum_{n=0}^{\infty}\widehat{\mathcal{G}}_{n}(v,w,z; \lambda, u) \frac{t^n}{n!}=\frac{(1-u)t}{\lambda e^{t}-u} e^{vt+wt^2+zt^3},
		\end{equation*}
		and  
		\begin{equation*}
			\sum_{n=0}^{\infty}\widehat{\mathcal{G}}_{n}(x,y,z; u) \frac{t^n}{n!}=\frac{(1-u)t}{e^{t}-u} e^{xt+yt^2+zt^3}.
		\end{equation*}
	\end{enumerate}
		
	\begin{thm}\label{thm-3-1} The multivariable generalized Hermite-type-Genocchi polynomials of order $a$ satisfy the relation
		\begin{align}
			&\widehat{\mathcal{G}}_{n,p}^{(k,a)}(x_1, x_2,\ldots, x_p; \nu, \lambda, u, c, d)\nonumber\\
			&\qquad=
			\sum_{m=0}^{n}\,\sum_{n_1+n_2+\ldots+n_p= m}(n)_m\widehat{\mathcal{G}}_{n-m}^{(k,a)}(\lambda, \nu, u, a, b)\,\prod_{j=1}^{p}\nu^{\frac{n_j}{j}} \binom{\frac{x_j}{\nu}}{\frac{n_j}{j}}\label{eq33new-01}
		\end{align} 
		Furthermore,
		\begin{align}
			&\widehat{\mathcal{G}}_{n,p}^{(k,a)}(x_1, x_2,\ldots, x_p; \nu, \lambda, u, c, d)\nonumber\\
			&\qquad=\sum_{s=0}^{\infty}\,\sum_{r_1+r_2+\ldots+r_p=s}\widehat{\mathcal{G}}_{n,p,\widetilde{w}}^{(k,a)}(n_1,n_2,\ldots,n_p;\nu, \lambda, u, c, d)\,{x_1}^{r_1}{x_2}^{r_2}\ldots  {x_p}^{r_p}.\label{eq33new-02}
		\end{align} 
		where
		\begin{align*}
			&\widehat{\mathcal{G}}_{n,p,\widetilde{w}}^{(k,a)}(n_1,n_2,\ldots,n_p;\nu, \lambda, u, c, d)\\
			&\qquad=\sum_{m=0}^{n}\,(n)_m\widehat{\mathcal{G}}_{n-m}^{(k,a)}(\nu, \lambda, u, c, d)\, \sum_{n_1+n_2+\ldots+n_p=m}\prod_{j=1}^{p}\frac{\widetilde{w}_{\nu}(n_j/j,r_j)}{\left(\frac{n_j}{j}\right)!}.
		\end{align*}
	\end{thm}
	
	\begin{proof} 
		By employing equation \eqref{eq13-1}, equation \eqref{eq711} givess
		\begin{align*}
			&\sum_{n=0}^{\infty}\widehat{\mathcal{G}}_{n,p}^{(k,a)}(x_1, x_2,\ldots, x_p; \nu, \lambda, u, c, d)\frac{t^n}{n!}\\
			&\;\;\;\;=\sum_{n=0}^{\infty}\left(\sum_{m=0}^{n}\binom{n}{m}\widehat{\mathcal{G}}_{n-m}^{(k,a)}(\nu, \lambda, u, c, d)\left\{\sum_{n_1+n_2+\ldots+n_p= m}\,m!\prod_{j=1}^{p}\nu^{\frac{n_j}{j}} \binom{\frac{x_j}{\nu}}{\frac{n_j}{j}}\right\}
			\right)\frac{t^n}{n!},
		\end{align*}
		which consequently yields equation \eqref{eq33new-01}. Moreover, using a certain variation of Stirling numbers of the first kind
		$$(x)_{n,m}=\sum_{j=0}^n\widetilde{w}_{m}(n,j)x^j,$$
		(see \cite{Mezo}), equation \eqref{eq33new-01} gives
			\begin{align*}
			&\sum_{n=0}^{\infty}\widehat{\mathcal{G}}_{n,p}^{(k,a)}(x_1, x_2,\ldots, x_p; \nu, \lambda, u, c, d)\frac{t^n}{n!}\\
			&\;\;\;\;=\sum_{n=0}^{\infty}\left(\sum_{m=0}^{n}\,\sum_{n_1+n_2+\ldots+n_p= m}(n)_m\widehat{\mathcal{G}}_{n-m}^{(k,a)}(\nu, \lambda, u, c, d) \prod_{j=1}^{p}\frac{1}{\left(\frac{n_j}{j}\right)!}\right.\\
			&\qquad\qquad\left.\times\,\sum_{s=0}^{\infty}\left\{\sum_{r_1,r_2,\ldots,r_p=s}\prod_{j=1}^{p}\widetilde{w}_{\nu}(n_j/j,r_j)\right\}{x_1}^{r_1}{x_2}^{r_2}\ldots  {x_p}^{r_p}\right)\frac{t^n}{n!}.
		\end{align*}
		Comparing the coefficients of $\frac{t^n}{n!}$ yields equation \eqref{eq33new-02}. 
	\end{proof}
	
	\begin{thm}\label{thm-3-2} The multivariable generalized Hermite-type-Genocchi polynomials of order $a$ satisfy the relation
		\begin{align*}
			&\widehat{\mathcal{G}}_{n}^{(k,a)}(x_1+x_2,y_1+y_2; \nu, \lambda, u, c, d)\\ 
			&\qquad=\sum_{q=0}^{n}\sum_{n_1+n_2+\ldots+n_p= n-q}\frac{(n-q)!\widehat{\mathcal{G}}_{q,p}^{(k,a)}(x_1,x_2,\ldots, x_p; \nu, \lambda, u, c, d)}{\prod_{j=1}^{p}{\left(\frac{n_j}{j}\right)!}}\,\,\left\{\prod_{j=1}^{p}(x_j){_{\frac{n_j}{j},\nu}}\right\}. 
		\end{align*}
	\end{thm}
	
	\begin{proof} 
		Employing equation \eqref{eq711} yields
		\begin{align*}
			&\sum_{n=0}^{\infty}\widehat{\mathcal{G}}_{n,p}^{(k,a)}(x_1+y_1, x_2+y_2,\ldots, x_p+y_p; \nu, \lambda, u, c, d)\frac{t^n}{n!}\\
			&=\left(\sum_{n=0}^{\infty}\widehat{\mathcal{G}}_{n,p}^{(k,a)}(x_1,x_2,\ldots, x_p; \nu, \lambda, u, c, d) \frac{t^n}{n!}\right)\left(\sum_{n\geq0}\left(\sum_{n_1+n_2+\ldots+n_p= n}\,n!\prod_{j=1}^{p}\nu^{\frac{n_j}{j}} \binom{\frac{x_j}{\nu}}{\frac{n_j}{j}}\right)\frac{t^n}{n!}\right)\\
			&=\sum_{n=0}^{\infty}\left(\sum_{q=0}^{n}(n-q)!\widehat{\mathcal{G}}_{q,p}^{(k,a)}(x_1,x_2,\ldots, x_p; \nu, \lambda, u, c, d) \sum_{n_1+n_2+\ldots+n_p= n-q}\,\prod_{j=1}^{p}\frac{(x_j){_{\frac{n_j}{j},\nu}}}{\left(\frac{n_j}{j}\right)!}\right)\frac{t^n}{n!}.
		\end{align*}
		Consequently, we obtain the desired result.
	\end{proof}

	\section{Acknowledgement}
	
	The authors sincerely acknowledge Cebu Normal University (CNU) for its financial support of this work through the Research Institute for Computational Mathematics and Physics (RICMP).

\end{document}